\documentclass[a4paper,11pt,twoside,leqno]{article}
\usepackage{amsmath,amssymb,amsfonts,amsthm,graphicx,fancyhdr}
\usepackage[latin1]{inputenc}
\usepackage[T1]{fontenc}
\usepackage{rotating,epic}
\usepackage[colorlinks=true]{hyperref}
\usepackage{enumitem}

\newtheorem{theorem}{Theorem}

\newtheorem{e-definition}[theorem]{Definition}

 \newtheoremstyle{drem}
      {3pt}
      {3pt}
      {\rmfamily}
      {}
      {\bf}
      {:}
      {.5em}
      {}

 \theoremstyle{drem}

\newcommand{\maxx}[1]{\textrm{\raisebox{.5ex}{\mbox{$\underset{#1}{\max}$}}} \:}
\newcommand{\zz}[0]{\mathbb{Z}}
\newcommand{\ie}[0]{\emph{i.e.}}

\newcommand{\comp}[0]{\mathsf{c}}
\newcommand{\bni}[0]{B_n^{\comp,\infty}}

\title{A remark on the connectedness of spheres in Cayley graphs.}
\date{\today}
\author{Antoine Gournay}

\begin{document}

\maketitle

\section{Introduction}\label{intro}

The aim of this note is to prove an elementary yet useful property of finitely presented groups. This property is called ``connected spheres'' in Blach\`ere's work \cite{Blachere} (where he shows that the Heisenberg group has this property). Filimonov \& Kleptsyn \cite{FK} use this remark to get some nice results on certain groups of diffeomorphisms of the circle.

Recall that, for a finitely generated group $\Gamma$ and $S \subset \Gamma$ a finite set such that $s \in S \implies s^{-1} \in S$, the \emph{Cayley graph} is the graph whose vertices are the elements of $\Gamma$ and where $g,h \in G$ are connected by an edge whenever there exists $s \in S$ such that $gs = h$. This $1$-complex is central to the study of $\Gamma$ as a geometric object.

A very rough property of Cayley graphs is the number of ends. Let $B_n$ be the ball of radius $n$ with centre at the identity element. This is defined to be the number of infinite connected components in the complement of $B_n$ as $n \to \infty$. Hopf \cite{Hopf} showed that a Cayley graph may have only $0$ (finite group), $1$, $2$ , or $\infty$ many ends. Stallings \cite{Stallings} described the case of groups with $2$ ends (virtually-$\zz$) and $\infty$ many ends (certain amalgamated products and HNN-extensions). Thus, it turns out ``most'' groups have $1$ end. 

The subject matter here is the number of ``important'' connected components in the spheres of thickness $r$. The term ``important'' needs to be added because the complement of $B_n$ may have many finite connected components (and only the infinite one is of interest here). The aim is to show that when the group is finitely presented, there exists $r$ (independent of $n$) such that these spheres are always connected. The complement of a set $A$ will be denoted $A^\comp$.
\begin{e-definition}\label{ladef}
Assume $\Gamma$ is one-ended (and finitely generated). Let $\bni$ be the the infinite connected component of $B_n^\comp$. For $r >0$, a graph has the property of connected spheres with constant $r$ if, for all $n \geq 0$, $B_{n+r} \cap \bni$ is connected. 
\end{e-definition}
When the constant is not specified, it should be interpreted that this is true for some $r$. It is necessary to restrict to the infinite connected component of $B_n^\comp$ because of dead-ends. See section \S{}\ref{sdead} below for further discussion on this topic.

Denote by $|w|$ the word length of a relation. 
\begin{theorem}\label{leteo}
Let $\Gamma$ be a finitely generated group with one end. Assume $\Gamma$ is finitely presented: $\Gamma= \langle S \mid R \rangle$. Take $r > \maxx{w \in R} \tfrac{|w|}{2}$. Then the Cayley graph of $\Gamma$ (with respect to generating set $S$) has connected spheres with constant $r$.
\end{theorem}

\par For completion, one could say that a non-empty subset $\Omega$ in a graph is simply connected if both $\Omega$ and its complement are connected. Let $\Omega^{+r}$ denote the set obtained by adding to $\Omega$ all points at distance $\leq r$ from $\Omega$. Then the above proof also carries to the following situation: in the Cayley graph of a finitely presented group, if $r > \tfrac{1}{2} \maxx{w \in R} |w|$ and $\Omega$ is simply connected, then $\Omega_n^{+r} \setminus \Omega_n$ is connected.

\emph{Addendum:} R.~Lyons pointed out to the author that the above theorem also follows from results of Babson \& Benjamini \cite{BB}. See also Tim{\'a}r \cite[Theorem 5.1]{Timar} and the book by Lyons (with Peres) \cite[Lemma 7.28]{Lyons}.

\emph{Acknowledgments:} The author would like to thanks V.~Kleptsyn for his encouragement to publish this small note and for pointing out its use in the paper of Filimonov \& Kleptsyn \cite{FK}. Comments from and discussions with J.~Brieussel, E.~Fink and J.~Lehnert helped improve this note. 

\emph{Remark:} The property of connected spheres was called ``uniformly one-ended'' in \cite[\S{}4.3]{moi-cohomlp-old}. This result was removed from subsequent versions of the paper since there was a mistake in its application, and the author could not find any interesting application. It then became clear from subsequent discussions with various people and from its use in the paper of Filimonov \& Kleptsyn \cite{FK} that, notwithstanding its elementary proof, this result is actually quite useful.

\section{The Cayley $2$-complex}\label{cayley}

When a group is finitely presented, one can associate the so-called Cayley $2$-complex $M_\Gamma$ to it. 
See Bridson \& Haefliger \cite[\S{}I.8A]{BH} for details.
Let $R$ be a (finite) set of (cyclically and ... reduced) relations associated to the (finite) generating set $S$. This complex is constructed as follows. Partition $S$ in sets of the form $A_i = \{s\} \cup \{s^{-1}\}$ where $i = 1, \ldots, n$. The $0$-skeleton is made of a single point $\star$. The $1$-skeleton is made of $n$ loops (with both ends at $\star$). Each of these loops is given an orientation and a label $a_i \in A_i$. This yields a bouquet of circles.

For each word $w = s_1s_2 \ldots s_k$ in $R$, take a disc whose boundary circle is cut into $k$ segments. The $j^\text{th}$ segment (in clockwise order) being labelled by the $a_i$ in $\{s_j\} \cup \{s_j^{-1}\}$ and oriented clockwise if $a_i = s_j$ and counter-clockwise otherwise. These discs are then glued, respecting orientation and label, to the bouquet of circles. 

In fact a group is the fundamental group of a CW-complex with finite $k$-skeletons for $k \leq 2$ if and only if it is finitely presented. In other words, it may always be assumed that the complex has no $k$-cells for $k >2$. This can be shown using the cellular approximation theorem. 

Another important remark is that a group generated by a symmetric finite set $S$ which has a uniform bound on the length of its relations is finitely presented. Indeed, if all relations are of length $\leq \ell$, then there are at most $|S|^\ell$ non-trivial reduced words with letters in $S$ of length $\leq \ell$.

\section{Proof}\label{proof}

Since $\Gamma$ is finitely presented, it is the fundamental group of its Cayley $2$-complex $M_\Gamma$. The $1$-skeleton of its universal covering, $\widetilde{M_\Gamma}$, is the Cayley graph of $\Gamma$. Take $r > \tfrac{1}{2} \maxx{w \in R} |w|$. Given two points $g$ and $g'$ of $B_{n+r} \cap \bni$, they can be joined by a path inside $B_{n+r}$ passing through the identity (since balls are connected) which is geodesic between $e$ and $g$ and between $e$ and $g'$. They can also be joined by a path $\gamma$ lying outside $B_n$.

Since $\widetilde{M_\Gamma}$ is simply connected, the loop obtained from these two paths may be filled in with a [combinatorial] disc $D$ of minimal [combinatorial] area. The boundary of $D$ is a relation $w$ (in bold lines above). Its decomposition into smaller discs corresponding to the $2$-cells (\ie~the defining relations) is the van Kampen diagram for $w$ (see \emph{e.g.} Bridson \& Haefliger \cite{BH} as above or Bridson \cite[Theorem 4.2.2]{Bridson}). 

Note a $2$-cell may not have a boundary $0$-cell both in $B_n$ and in $B_{n+r}^\comp$. Indeed, this would imply that the length of its boundary word is $\geq 2r$ (since any path from $B_n$ to the complement of $B_{n+r}$ is of length at least $r$) and would contradict the choice of $r$: $2r > \maxx{w \in R} |w|$. 

\begin{center}
\includegraphics[width=10cm]{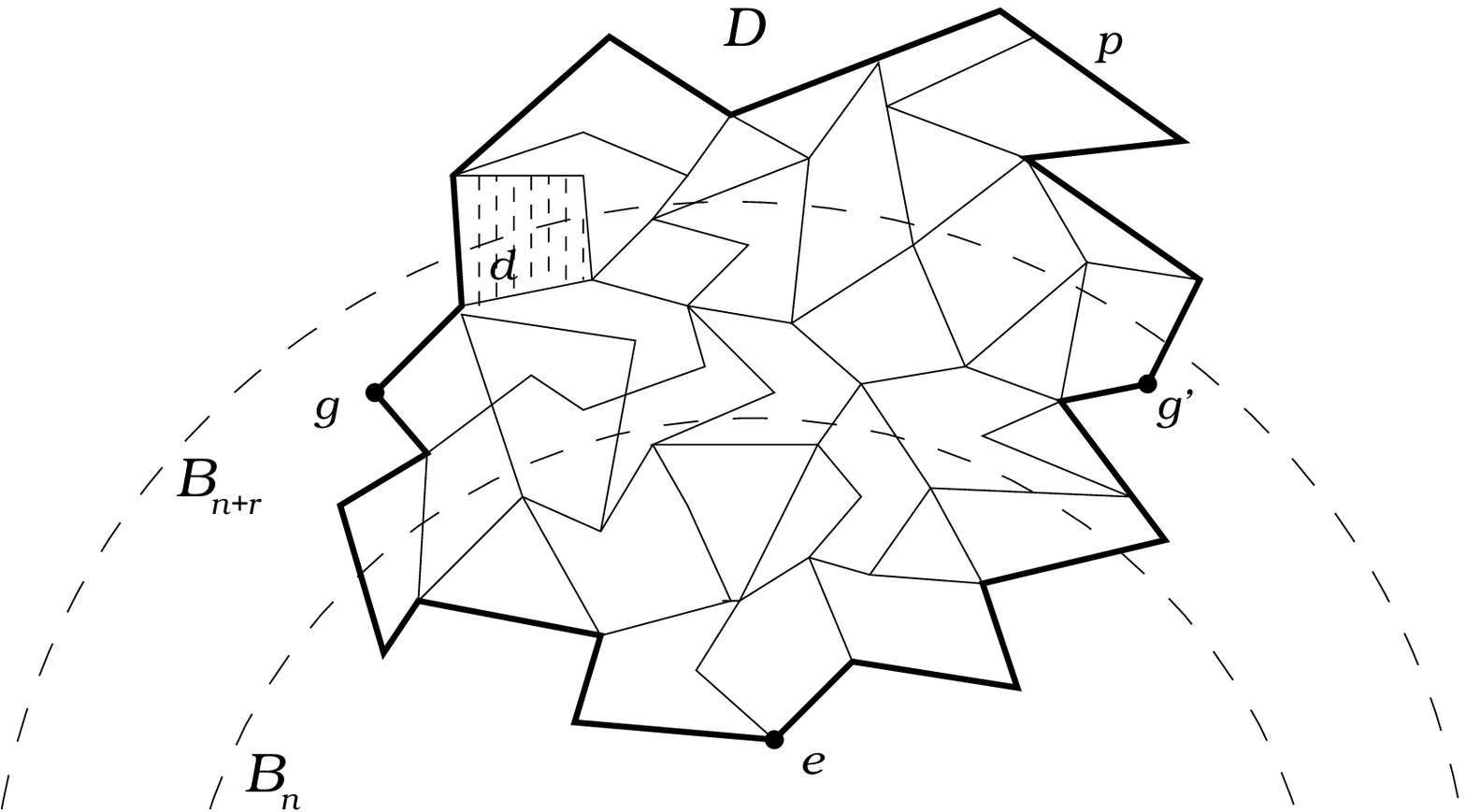}  
\end{center}

Take $p$ so that $D$ has minimal [combinatorial] area. Assume there is a boundary $2$-cell $d$ of $D$ which contains a $0$-cell in $B_{n+r}^\comp$. Upon removing $d$, $D$ might become disconnected. If this is the case, consider $D'$ (still a [combinatorial] disc) the connected component of $D \setminus d$ containing $e$. Its boundary may be used to define a new path $p'$ which contradicts the minimality of $p$. Indeed, by the the previous paragraph, $p'$ still lies outside $B_n$. 

Thus it may be assumed that $p$ does not contain any $0$-cell outside $B_{n+r}$ and lies outside $B_n$. This proves the claim.
\hfill \emph{Q.E.D}

\section{Dead-end?... Questions and further comments}\label{sdead}

What really matters for the connectedness of spheres is the retreat depth (or strong depth) of $\gamma \in \Gamma$ (for a generating set $S$). This is the smallest $d$ such that $\gamma$ is in $B^{\comp,\infty}_{|\gamma| -d}$ where $|\gamma|$ is the word length of $\gamma$. Lehnert \cite{Lehnert} (where it bears the name ``strong depth'') shows that for the Houghton group $H_2$ (a group which is not $FP_2$, hence not finitely presented) it is unbounded. Warshall \cite[Proposition 2]{Warshall} (where it bears the name ``retreat depth'') shows it is bounded for the Heisenberg group (for any generating set). 

J.~Lehnert pointed out to the author that retreat depth is not invariant under changing the (finite) generating set. The counterexample comes from lamplighter groups. Define the [usual] depth of an element $g$ to be the distance between $g$ and $B_{|g|}^{\comp,\infty}$. In \cite{Warshall2}, Warshall shows there is a generating set $S$ for which the lamplighter (on $\zz$) has bounded [usual] depth, hence bounded retreat depth. On the other hand, Cleary and Taback \cite{CT} describe dead-end elements (for the usual generators) which are readily seen to be of unbounded retreat depth.

A discussion with J.~Brieussel made it quite obvious that the lamplighter on $\zz$ (\ie $\zz_2 \wr \zz$) does not have connected spheres. This is no longer so obvious on $\zz_2 \wr \zz^2$. Funnily, $\zz \wr \zz$ (which does not have dead-ends with the usual generating set) has connected spheres. 

Here are a few interesting questions (which we believe should not be hard to prove or disprove). A group has $F_n$ if its $K(\Gamma,1)$ is finite in dimensions $\leq n$. Finitely presented is equivalent to $F_2$. Recall that a group has $FP_n$ (for a ring $R$) if there is a [partial] projective resolution of length $n$ by finitely generated $R\Gamma$-modules of the ring $R$. Finite presentation implies $FP_2$, but the converse is [non-trivially] false. It is usually understood that $R= \zz$, but in the following questions, it is not clear if a specific ring should be taken.
\begin{enumerate}[label=(\roman{enumi}), ref=(\roman{enumi})]
\item Does $FP_2$ implies connected spheres?
\item Is uniformly bounded retreat depth invariant of the generating set amongst groups with a finite presentation? \label{retde}
\item Is connected spheres invariant under changing the generating set? \label{fink}
\item If $\Gamma$ is such that $K(\Gamma,1)$ is finite, is the retreat depth uniformly bounded? \label{kg1fini}
\item Can one relax ``finite $K(\Gamma,1)$'' to $F_k$ or $FP_k$ (for some $k$) in \ref{kg1fini}? \label{finpres}
\item For a finitely generated group $\Gamma$, does there exist a $\alpha \in \{0,1,2,\infty\}$ and $r_0(S) >0$ ($r_0$ depends on the generating set), such that for any $r\geq r_0(S)$ the number of connected components of $B_{n+r} \cap \bni$ tends to $\alpha$ as $n \to \infty$?
\end{enumerate}
J.~Lehnert pointed out to the author that realistic candidates for a negative answer to question \ref{retde} and \ref{finpres} are Houghton's groups ($H_k$ is finitely presented for $k\geq 3$, has $FP_{k-1}$ but not $FP_k$). \ref{fink} was pointed out to the author by E.~Fink.

Lastly, it might be interesting to (try to) generalise the above result to higher filling properties and groups with property $F_n$ or $FP_n$.





\end{document}